\documentclass[12pt,reqno]{amsart}

\usepackage{amssymb}

\usepackage{eucal}

\input{diagrams}

\font\sss=cmss8

\def\cA{{\mathcal A}}
\def\cB{{\mathcal B}}

\def\sD{\mbox{\sf D}}

\def\D{\sD}

\def\depth{\operatorname{depth}}

\def\dim{\operatorname{dim}}
\def\Dsmall{\mbox{\sss D}}

\def\Ext{\operatorname{Ext}}

\def\fd{\operatorname{fd}}

\def\H{\operatorname{H}}

\def\Hom{\operatorname{Hom}}
\def\id{\operatorname{id}}

\def\inf{\operatorname{inf}}

\def\limdir{\operatorname{lim}_{_{_{_{\hspace{-17pt} 
                                      \textstyle \longrightarrow}}}}}

\def\LTensor{{\otimes}^{\operatorname{L}}}

\def\pd{\operatorname{pd}}

\def\RHom{\operatorname{RHom}}

\def\sup{\operatorname{sup}}

\def\width{\operatorname{width}}
\def\Z{\operatorname{Z}}

\numberwithin{equation}{part}

\newtheorem{Lemma}{Lemma}[section]
\newtheorem{Theorem}[Lemma]{Theorem}
\newtheorem{Proposition}[Lemma]{Proposition}

\theoremstyle{definition}
\newtheorem{Definition}[Lemma]{Definition}
\newtheorem{Setup}[Lemma]{Setup}

\newtheorem{Remark}[Lemma]{Remark}

\def\skewtimes{\ltimes}

\def\Gpd{\operatorname{Gpd}}
\def\Gfd{\operatorname{Gfd}}
\def\Gid{\operatorname{Gid}}

\def\CGpd{\Gpd_{A \skewtimes C}}
\def\CGfd{\Gfd_{A \skewtimes C}}
\def\CGid{\Gid_{A \skewtimes C}}

\def\Gdim{\operatorname{Gdim}}
\def\CGdim{\mbox{$C$-$\Gdim$}}

\def\CMpd{\operatorname{CMpd}}
\def\CMfd{\operatorname{CMfd}}
\def\CMid{\operatorname{CMid}}

\def\CMdim{\operatorname{CMdim}}

\begin{document}

\title[Cohen-Macaulay dimensions]
{Cohen-Macaulay injective, projective, and flat dimension}

\author{Henrik Holm \ \ }
\address{Matematisk Afdeling, Universitetsparken 5, 2100 K\o benhavn
\O, Denmark}
\email{holm@math.ku.dk}

\author{\ \ Peter J\o rgensen}
\address{Department of Pure Mathematics, University of Leeds,
Leeds LS2 9JT, United Kingdom}
\email{popjoerg@maths.leeds.ac.uk, www.maths.leeds.ac.uk/\~{ }popjoerg}

\keywords{Semi-dualizing module, Gorenstein homological dimension,
Cohen-Macaulay homological dimension}

\subjclass[2000]{13D05, 13D25}

\begin{abstract} 

We define three new homological dimensions --- Co\-hen-Ma\-cau\-lay
injective, projective, and flat dimension --- which inhabit a theory
similar to that of classical injective, projective, and flat
dimension.  Finiteness of the new dimensions characterizes
Cohen-Macaulay rings with dualizing modules.

\end{abstract}

\maketitle

\setcounter{section}{0}
\section{Introduction}
\label{sec:intro}

The classical theory of injective, projective, and flat dimension has
had great success in commutative algebra.  In particular, it has been
very useful to know that finiteness of these dimensions characterizes
regular rings.

Several attempts have been made to mimic this success by constructing
homological dimensions whose finiteness would characterize other rings
than the regular ones.  These efforts have given us complete
intersection dimension, Gorenstein dimension, and Cohen-Macaulay
dimension.

The normal practice has not been to mimic all three classical
dimensions, but rather to focus on projective dimension for finitely
generated modules.  Hence complete intersection dimension and
Cohen-Macaulay dimension only exist in this restricted sense, and the
same used to be the case for Gorenstein dimension.

However, recent years have seen much work on the Gorenstein theory
which now contains both Gorenstein injective, projective, and flat
dimension.  These dimensions inhabit a nice theory similar to the
classical one.  A good summary is in \cite{Wintherbook}, although this
reference is already a bit dated.

The purpose of this paper is to do something similar in the
Cohen-Macaulay case.  So we define Cohen-Macaulay injective,
projective, and flat dimension, and show some central properties.  

Our main result is theorem \ref{thm:CMdims_finite} which lists a large
number of finiteness conditions on the Cohen-Macaulay dimensions, and
shows that they are all equivalent to the ground ring being
Cohen-Macaulay with a dualizing module.  As a sample of further
possibilities, there is also an Auslander-Buchsbaum formula for
Cohen-Macaulay projective dimension, see the\-o\-rem \ref{thm:AB}.

As tools to define the Cohen-Macaulay dimensions, we use ``ring
changed'' Gorenstein homological dimensions.  If $A$ is a ring with a
semi-dualizing module $C$ (as defined in \cite{Wintherpaper}), then we
can consider the trivial extension ring $A \skewtimes C$, and if $M$
is a complex of $A$-modules, then we can consider $M$ as a complex of
$(A \skewtimes C)$-modules and take ``ring changed'' Gorenstein
dimensions of $M$ over $A \skewtimes C$.  We shall develop the theory
of these dimensions further in \cite{HHPJ1}.  The present paper only
refers to results from \cite{HHPJ1} at the end, in the proofs of lemma
\ref{lem:Gerko_inequality} and theorem \ref{thm:CMpd_and_CMdim}.

The paper is organized as follows: Section \ref{sec:CGdims} defines
the Cohen-Macaulay dimensions.  Section \ref{sec:extension} studies
the trivial extension ring $A \skewtimes C$ when $C$ is a
semi-dualizing module for $A$.  Section \ref{sec:estimates} gives some
bounds on the injective dimension of $C$.  And finally, section
\ref{sec:CMdims} studies the Cohen-Macaulay dimensions and proves the
results we have stated.

\begin{Setup}
Throughout, $A$ is a commutative noetherian ring.  Complexes of
$A$-modules have the form
\[
  \cdots \longrightarrow M_{i+1}
  \longrightarrow M_i
  \longrightarrow M_{i-1}
  \longrightarrow \cdots,
\]
and the words ``right-bounded'' and ``left-bounded'' are to be
understood relative to this.
\end{Setup}

\section{Cohen-Macaulay dimensions}
\label{sec:CGdims}

\begin{Definition}
Let $C$ be an $A$-module.  The direct sum $A \oplus C$ can be equipped
with the product
\[
  \left( \begin{array}{c} a \\ c \end{array} \right) 
  \cdot \left( \begin{array}{c} a_1 \\ c_1 \end{array} \right)
  = \left( \begin{array}{c} aa_1 \\ ac_1 + ca_1 \end{array} \right) .
\]
This turns $A \oplus C$ into a ring which is called the trivial
extension of $A$ by $C$ and denoted $A \skewtimes C$.
\end{Definition}

There are ring homomorphisms
\[
  \begin{array}{ccccc}
    A & \longrightarrow & A \skewtimes C & \longrightarrow & A, \\[.3cm]
    a & \longmapsto
      & \left( \begin{array}{c} a \\ 0 \end{array} \right) \lefteqn{,}
      & 
      & \\[.6cm]
      &
      & \left( \begin{array}{c} a \\ c \end{array} \right) 
      & \longmapsto 
      & a
  \end{array}
\] 
whose composition is the identity on $A$.  These homomorphisms allow
us to view any $A$-module as an $(A \skewtimes C)$-module and any $(A
\skewtimes C)$-module as an $A$-module, and we shall do so freely.

In particular, if $M$ is a complex of $A$-modules with suitably
bounded homology, then we can consider the ``ring changed'' Gorenstein
homological dimensions
\[
  \CGid M, \; \CGpd M, \; \mbox{and} \; \CGfd M,
\]
where $\Gid$, $\Gpd$, and $\Gfd$ denote the Gorenstein injective,
projective, and flat dimensions, as described e.g.\ in
\cite{Wintherbook}. 

Before the next definition, recall that a semi-dualizing module $C$
for $A$ is a finitely generated module for which the canonical map $A
\longrightarrow \Hom_A(C,C)$ is an isomorphism, while $\Ext_A^i(C,C) =
0$ for each $i \geqslant 1$.  Equivalently, $C$ is a finitely generated
module so that the canonical morphism $A \longrightarrow \RHom_A(C,C)$
is an isomorphism in the derived category $\D(A)$.  An example of a
semi-dualizing module is $A$ itself.  The theory of semi-dualizing
modules (and complexes) is developed in \cite{Wintherpaper}.

\begin{Definition}
\label{def:CMdims}
Let $M$ and $N$ be complexes of $A$-modules so that the homology of
$M$ is bounded to the left and the homology of $N$ is bounded to the
right.

The Cohen-Macaulay injective, projective, and flat dimensions of $M$
and $N$ over $A$ are
\begin{align*}
  \CMid_A M & = \inf \, \{\, \CGid M \;|\; C \mbox{ is a semi-dualizing module} \,\}, \\
  \CMpd_A N & = \inf \, \{\, \CGpd N \;|\; C \mbox{ is a semi-dualizing module} \,\}, \\
  \CMfd_A N & = \inf \, \{\, \CGfd N \;|\; C \mbox{ is a semi-dualizing module} \,\}. \\
\end{align*}
\end{Definition}

\section{Lemmas on the trivial extension}
\label{sec:extension}

\begin{Lemma}
\label{lem:induced_injectives}
Let $C$ be an $A$-module.
\begin{enumerate}

  \item  If $I$ is a (faithfully) injective $A$-module then $\Hom_A(A
         \skewtimes C,I)$ is a (faithfully) injective $(A \skewtimes
         C)$-module.

\smallskip

  \item  Each injective $(A \skewtimes C)$-module is a direct summand
         in a module  $\Hom_A(A \skewtimes C,I)$ where $I$ is an
         injective $A$-module.

\end{enumerate}
\end{Lemma}

\begin{proof}
(1)  Adjunction gives
\begin{align}
\nonumber
  \Hom_{A \skewtimes C}(-,\Hom_A(A \skewtimes C,I))
  & \simeq \Hom_A((A \skewtimes C) \otimes_{A \skewtimes C} -,I) \\
\label{equ:adjunction}
  & \simeq \Hom_A(-,I)
\end{align}
on $(A \skewtimes C)$-modules, making it clear that if $I$ is a
(faithfully) injective $A$-module, then $\Hom_A(A \skewtimes C,I)$ is a
(faithfully) injective $(A \skewtimes C)$-module.

\medskip
\noindent
(2)  To see that an injective $(A \skewtimes C)$-module $J$ is a
direct summand in a module of the form $\Hom_A(A \skewtimes C,I)$, it is
enough to embed it into such a module.  For this, first view $J$ as an
$A$-module and embed it into an injective $A$-module $I$.  Then use
equation \eqref{equ:adjunction} to convert the monomorphism of
$A$-modules $J \hookrightarrow I$ to a monomorphism of $(A \skewtimes
C)$-modules $J \hookrightarrow \Hom_A(A \skewtimes C,I)$.
\end{proof}

\begin{Lemma}
\label{lem:extension_formulae}
Let $C$ be a semi-dualizing module for $A$.
\begin{enumerate}

  \item  There is an isomorphism 
\[
  \RHom_A(A \skewtimes C,C) \cong A \skewtimes C 
\]
         in the derived category $\D(A \skewtimes C)$.

\smallskip

  \item  There is a natural equivalence 
\[
  \RHom_{A \skewtimes C}(-,A \skewtimes C) \simeq \RHom_A(-,C)
\]
         of functors on $\D(A)$.

\smallskip

  \item  If $M$ is in $\D(A)$ then the biduality morphisms  
\[
  M \longrightarrow \RHom_A(\RHom_A(M,C),C)
\]
         and
\[
  M \longrightarrow \RHom_{A \skewtimes C}(
    \RHom_{A \skewtimes C}(M,A \skewtimes C),A \skewtimes C
                                          )
\]
         are equal.

\smallskip

  \item  There is an isomorphism 
\[
  \RHom_{A \skewtimes C}(A,A \skewtimes C) \cong C
\]
         in $\D(A \skewtimes C)$.

\end{enumerate}
\end{Lemma}

\begin{proof}
(1)  Since $C$ is semi-dualizing, it is clear that there is an
isomorphism in $\D(A)$,
\[
  \RHom_A(A \oplus C,C) \cong C \oplus A.
\]
It is easy to see that in fact, this isomorphism respects the action
of $A \skewtimes C$, so
\[
  \RHom_A(A \skewtimes C,C) \cong A \skewtimes C
\]
in $\D(A \skewtimes C)$.

\medskip
\noindent
(2)  This is a computation,
\begin{align*}
  \RHom_{A \skewtimes C}(-,A \skewtimes C)
  & \stackrel{\rm (a)}{\simeq}
    \RHom_{A \skewtimes C}(-,\RHom_A(A \skewtimes C,C)) \\
  & \stackrel{\rm (b)}{\simeq}
    \RHom_A((A \skewtimes C) \LTensor_{A \skewtimes C} -,C) \\
  & \simeq \RHom_A(-,C),
\end{align*}
where (a) is by part (1) and (b) is by adjunction.

\medskip
\noindent
(3) and (4)  These are easy to obtain from (2).
\end{proof}

\begin{Lemma}
\label{lem:A_and_C_Gproj}
Let $C$ be a semi-dualizing module for $A$ and let $I$ be an injective
$A$-module.
\begin{enumerate}

  \item  $A$ and $C$ are Gorenstein projective over $A \skewtimes C$.

\smallskip

  \item  $\Hom_A(A,I) \cong I$ and $\Hom_A(C,I)$ are Gorenstein
         injective over $A \skewtimes C$.
\end{enumerate}
\end{Lemma}

\begin{proof}
(1) Lemma \ref{lem:extension_formulae}(4) says
$\RHom_{A \skewtimes C}(A,A \skewtimes C) \cong C$.  That is, the dual
of $A$ with respect to the ring $A \skewtimes C$ is $C$.  But
dualization with respect to the ring preserves the class of finitely
generated Gorenstein projective modules by \cite[obs.\
(1.1.7)]{Wintherbook}, so to prove part (1) it is enough to see that
$A$ is Gorenstein projective over $A \skewtimes C$.

By \cite[prop.\ (2.2.2)]{Wintherbook}, this will follow if $\RHom_{A
\skewtimes C}(A,A \skewtimes C)$ is concentrated in degree zero and
the biduality morphism
\[
  A \longrightarrow 
  \RHom_{A \skewtimes C}(
    \RHom_{A \skewtimes C}(A,A \skewtimes C),A \skewtimes C
                        )
\]
is an isomorphism.  The first of these conditions holds by lemma
\ref{lem:extension_formulae}(4), and the second condition holds
because the biduality morphism equals
\[
  A \longrightarrow \RHom_A(\RHom_A(A,C),C)
\]
by lemma \ref{lem:extension_formulae}(3), and this is an isomorphism
because it is equal to the canonical morphism $A \longrightarrow
\RHom_A(C,C)$. 

\medskip
\noindent
(2) We will prove that $\Hom_A(C,I)$ is Gorenstein
injective over $A \skewtimes C$, the case of $\Hom_A(A,I) \cong I$
being similar.

Since $C$ is Gorenstein projective over $A \skewtimes C$, by
definition it has a complete projective resolution $P$.  So $P$ is a
complex of $(A \skewtimes C)$-modules which has $C$ as one of its
cycle modules, say $\Z_0(P) \cong C$.  Moreover, $P$ is an exact
complex of projective $(A \skewtimes C)$-modules, and
\[
  \Hom_{A \skewtimes C}(P,Q) 
\]
is exact when $Q$ is a projective $(A \skewtimes C)$-module.  Since
$C$ is finitely generated, we can assume that $P$ consists of finitely
generated $(A \skewtimes C)$-modules by \cite[thms.\ (4.1.4) and
(4.2.6)]{Wintherbook}.

The $(A \skewtimes C)$-module $J = \Hom_A(A \skewtimes C,I)$ is
injective by lemma \ref{lem:induced_injectives}(1).  Consider the
complex
\[
  K = \Hom_{A \skewtimes C}(P,J).
\]
This is clearly an exact complex of injective $(A \skewtimes
C)$-modules.  Moreover, if $L$ is an injective $(A
\skewtimes C)$-module then
\begin{align*}
  \Hom_{A \skewtimes C}(L,K)
  & = \Hom_{A \skewtimes C}(L,\Hom_{A \skewtimes C}(P,J)) \\
  & \cong \Hom_{A \skewtimes C}(L \otimes_{A \skewtimes C} P,J) \\
  & \cong \Hom_{A \skewtimes C}(P,\Hom_{A \skewtimes C}(L,J)) \\
  & = (*),
\end{align*}
where both $\cong$'s are by adjunction.  Here $\Hom_{A \skewtimes
C}(L,J)$ is a flat $(A \skewtimes C)$-module by \cite[thm.\
1.2]{Lazard}, so it is the direct limit of projective $(A \skewtimes
C)$-modules,
\[
  \Hom_{A \skewtimes C}(L,J) \cong \limdir Q_{\alpha}.
\]
So
\[
  (*) \cong \Hom_{A \skewtimes C}(P,\limdir Q_{\alpha})
      \cong \limdir \Hom_{A \skewtimes C}(P,Q_{\alpha}) = (**),
\]
where the second $\cong$ holds because each module in $P$ is finitely
generated.  Since each $\Hom_{A \skewtimes C}(P,Q_{\alpha})$ is exact,
so is $(**)$.

This shows that $K$ is a complete injective resolution over 
$A \skewtimes C$, and 
\begin{align*}
  \Z_{-1}(K) & = \Z_{-1}(\Hom_{A \skewtimes C}(P,J)) \\
  & \cong \Hom_{A \skewtimes C}(\Z_0(P),J) \\
  & \cong \Hom_{A \skewtimes C}(C,\Hom_A(A \skewtimes C,I)) \\
  & \stackrel{\rm (a)}{\cong} 
    \Hom_A(C \otimes_{A \skewtimes C} (A \skewtimes C),I) \\
  & \cong \Hom_A(C,I), \\
\end{align*}
where (a) is again by adjunction.  So $K$ is a complete
injective resolution of $\Hom_A(C,I)$ which is therefore Gorenstein
injective. 
\end{proof}

\begin{Lemma}
\label{lem:hom_from_injective}
Let $C$ be a semi-dualizing module for $A$ and let $I$ be an injective
$A$-module.  Then there is an equivalence of functors on $\D(A
\skewtimes C)$,
\[
  \RHom_{A \skewtimes C}(\Hom_A(A \skewtimes C,I),-)
    \simeq \RHom_A(\Hom_A(C,I),-).
\]
\end{Lemma}

\begin{proof}
This is a computation,
\begin{align*}
  \lefteqn{\RHom_{A \skewtimes C}(\Hom_A(A \skewtimes C,I),-)}
    \;\;\;\;\; & \\
  & \simeq \RHom_{A \skewtimes C}(\RHom_A(A \skewtimes C,I),-) \\
  & \stackrel{\rm (a)}{\simeq} 
    \RHom_{A \skewtimes C}(
    \RHom_A(\RHom_A(A \skewtimes C,C),I),-
                          ) \\
  & \stackrel{\rm (b)}{\simeq}
    \RHom_{A \skewtimes C}(
    (A \skewtimes C) \LTensor_A \RHom_A(C,I),-
                          ) \\
  & \stackrel{\rm (c)}{\simeq}
    \RHom_A(\RHom_A(C,I),-), \\
  & \simeq \RHom_A(\Hom_A(C,I),-) \\
\end{align*}
where (a) is by lemma \ref{lem:extension_formulae}(1), (b) is by
\cite[(A.4.24)]{Wintherbook}, and (c) is by adjunction.
\end{proof}

\begin{Lemma}
\label{lem:AxC_Gorenstein}
Assume that the ring $A$ is local and let $C$ be a finitely generated
$A$-module.  Then
\[
  \mbox{ $A \skewtimes C$ is a Gorenstein ring
         $\Leftrightarrow$
         $C$ is a dualizing module for $A$.}
\]
\end{Lemma}

\begin{proof}
This can be found between the lines in \cite{FoxbyGorenstein} or
\cite{ReitenGorenstein}, or explicitly as a special case of
\cite[thm.\ 2.2]{PJRecognizing}.
\end{proof}

\section{Bounds on the injective dimension of $C$}
\label{sec:estimates}

\begin{Lemma}
\label{dualt_PJ-Lemma10}
Let $C$ be a semi-dualizing module for $A$ and let $M$ be an
$A$-module which is Gorenstein injective over $A \skewtimes C$. Then
there exists a short exact sequence of $A$-modules,
\[
  0 \rightarrow M^{\prime} 
    \longrightarrow \Hom_A(C,I) 
    \longrightarrow M
    \rightarrow 0, 
\]
where $I$ is injective over $A$ and $M^{\prime}$ is Gorenstein
injective over $A \skewtimes C$, which stays exact if one applies to
it the functor $\Hom_A(\Hom_A(C,J),-)$ for any injective $A$-module
$J$.
\end{Lemma}

\begin{proof}
Since $M$ is Gorenstein injective over $A \skewtimes C$, it has a
complete injective resolution.  From this can be extracted a short
exact sequence of $(A \skewtimes C)$-modules,
\[
  0 \rightarrow N
    \longrightarrow K
    \longrightarrow M
    \rightarrow 0,
\]
where $K$ is injective and $N$ Gorenstein injective over $A
\skewtimes C$, which stays exact if one applies to it the functor
$\Hom_{A \skewtimes C}(L,-)$ for any injective $(A \skewtimes
C)$-module $L$.  

In particular, the sequence stays exact if one applies to it the
functor $\Hom_{A \skewtimes C}(\Hom_A(A \skewtimes C,J),-)$ for any
injective $A$-module $J$, because $\Hom_A(A \skewtimes C,J)$ is an
injective $(A \skewtimes C)$-module by lemma
\ref{lem:induced_injectives}(1).

By lemma \ref{lem:induced_injectives}(2), the injective $(A
\skewtimes C)$-module $K$ is a direct summand in $\Hom_A(A \skewtimes
C,I)$ for some injective $A$-module $I$.  If $K \oplus K^{\prime}
\cong \Hom_A(A \skewtimes C,I)$, then by adding $K^{\prime}$ to both
the first and the second module in the short exact sequence, we
may assume that the sequence has the form
\[
  0 \rightarrow N
    \longrightarrow \Hom_A(A \skewtimes C,I)
    \stackrel{\eta}{\longrightarrow} M
    \rightarrow 0.
\]
The module $N$ is still Gorenstein injective over $A \skewtimes C$,
and the sequence still stays exact if one applies to it the functor
\[
  \Hom_{A \skewtimes C}(\Hom_A(A \skewtimes C,J),-)
\]
for any injective $A$-module $J$.

Now let us consider in detail the homomorphism $\eta$.  Elements of
the source $\Hom_A(A \skewtimes C,I)$ have the form $(\alpha,\gamma)$
where $A \stackrel{\alpha}{\longrightarrow} I$ and $C
\stackrel{\gamma}{\longrightarrow} I$ are homomorphisms of
$A$-modules.  The $(A \skewtimes C)$-module structure of $\Hom_A(A
\skewtimes C,I)$ comes from the first variable, and one checks that it
takes the form
\[
  \left(
    \begin{array}{c}
      a \\ c
    \end{array}
  \right)
  \cdot
  (\alpha,\gamma)
  =
  (a \alpha + \chi_{\gamma(c)},a \gamma),
\]
where $\chi_{\gamma(c)}$ is the homomorphism $A \longrightarrow I$
given by $a \mapsto a\gamma(c)$.

The target of $\eta$ is $M$ which is an $A$-module.  When viewed as an
$(A \skewtimes C)$-module, $M$ is annihilated by the ideal $0
\skewtimes C$, so
\begin{equation}
\label{equ:eta1}
  0 
  =   
  \left(
    \begin{array}{c}
      0 \\ c
    \end{array}
  \right)
  \cdot
  \eta(\alpha,\gamma)
  =
  \eta(
  \left(
    \begin{array}{c}
      0 \\ c
    \end{array}
  \right)
  \cdot
  (\alpha,\gamma)
      )
  =
  \eta(\chi_{\gamma(c)},0),
\end{equation}
where the last $=$ follows from the previous equation.

In fact, this implies
\begin{equation}
\label{equ:eta2}
  \eta(\alpha,0) = 0
\end{equation}
for each $A \stackrel{\alpha}{\longrightarrow} I$.  To see so, note
that there is a surjection $F \longrightarrow \Hom_A(C,I)$ with $F$
free, and hence a surjection $C \otimes_A F \longrightarrow C
\otimes_A \Hom_A(C,I)$.  The target here is isomorphic to $I$ by
\cite[prop.\ (4.4) and obs.\ (4.10)]{Wintherpaper}, so there is a
surjection $C \otimes_A F \longrightarrow I$.  As $C \otimes_A F$ is a
direct sum of copies of $C$, this means that, given an element $i$ in
$I$, it is possible to find homomorphisms $\gamma_1, \ldots, \gamma_t
: C \longrightarrow I$ and elements $c_1, \ldots, c_t$ in $C$ with
$i = \gamma_1(c_1) + \cdots + \gamma_t(c_t)$.  Hence the homomorphism
\[
  A \stackrel{\alpha}{\longrightarrow} I 
\]
given by $a \mapsto ai$ is equal to
\[
  \chi_{\gamma_1(c_1) + \cdots + \gamma_t(c_t)}
  = \chi_{\gamma_1(c_1)} + \cdots + \chi_{\gamma_t(c_t)},
\]
and so equation \eqref{equ:eta1} implies equation \eqref{equ:eta2}. 

To make use of this, observe that the exact sequence of $(A \skewtimes
C)$-modules 
\begin{equation}
\label{equ:basic_exact_sequence}
  0 \rightarrow C \longrightarrow A \skewtimes C \longrightarrow A
  \rightarrow 0
\end{equation}
induces an exact sequence
\[
  0 \rightarrow \Hom_A(A,I)
    \longrightarrow \Hom_A(A \skewtimes C,I)
    \longrightarrow \Hom_A(C,I)
    \rightarrow 0.
\] 
So equation \eqref{equ:eta2} can be interpreted as saying that the
homomorphism $\Hom_A(A \skewtimes C,I)
\stackrel{\eta}{\longrightarrow} M$ factors through the surjection $
\Hom_A(A \skewtimes C,I) \longrightarrow
\Hom_A(C,I)$.  This means that we can construct a commutative diagram
of $(A \skewtimes C)$-modules with exact rows,
\[
  \begin{diagram}[labelstyle=\scriptstyle,midshaft]
    0 & \rTo & N & \rTo & \Hom_A(A \skewtimes C,I) & \rTo^{\eta} & M & \rTo & 0 \\
      & & \dTo & & \dTo & & \vEq & & \\
    0 & \rTo & M^{\prime} & \rTo & \Hom_A(C,I) & \rTo & M & \rTo & 0 \lefteqn{.} \\
  \end{diagram}
\]

We will show that if we view the lower row as a sequence of
$A$-modules, then it is a short exact sequence with the properties
claimed in the lemma.

First, $I$ is injective over $A$ by construction.

Secondly, applying the Snake Lemma to the above diagram embeds the
vertical arrows into exact sequences.  The leftmost of these is
\[
  0 \rightarrow \Hom_A(A,I) 
    \longrightarrow N
    \longrightarrow M^{\prime}
    \rightarrow 0.
\] 
Here the modules $\Hom_A(A,I) \cong I$ and $N$ are Gorenstein
injective over $A \skewtimes C$ by lemma \ref{lem:A_and_C_Gproj}(2),
respectively, by construction.  Hence $M^{\prime}$ is also Gorenstein
injective over $A \skewtimes C$ because the class of Gorenstein
injective modules is injectively resolving by \cite[thm.\
2.6]{HHGorensteinHomDim}.

Thirdly, by construction, the upper sequence in the diagram stays
exact if one applies to it the functor $\Hom_{A \skewtimes C}(\Hom_A(A
\skewtimes C,J),-)$ for any injective $A$-module $J$.  It follows
that the same holds for the lower row.  But taking $\H_0$ of
the isomorphism in lemma \ref{lem:hom_from_injective} shows 
\[
  \Hom_{A \skewtimes C}(\Hom_A(A \skewtimes C,J),-)
    \simeq \Hom_A(\Hom_A(C,J),-),
\]
so the lower row in the diagram also stays exact if one applies to it
the functor $\Hom_A(\Hom_A(C,J),-)$ for any injective $A$-module $J$.
\end{proof}

The following lemmas use ${}_{C}\cA(A)$ and ${}_{C}\cB(A)$, the
Auslander and Bass classes of the semi-dualizing module $C$, as
introduced in \cite[def.\ (4.1)]{Wintherpaper}.  The proof of the
first of the lemmas can be found in \cite{FHS}.

\begin{Lemma}
\label{lem:triangle}
Let $C$ be a semi-dualizing module for $A$, let $M$ in ${}_{C}\cA(A)$
satisfy $\CGid M < \infty$, and write $s = \sup \{\, i \,|\, \H_i\!M
\not= 0 \,\}$.  Then there is a distinguished triangle in $\D(A)$,
\[
  \Sigma^s H \longrightarrow Y \longrightarrow M \longrightarrow,
\]
where $H$ is an $A$-module which is Gorenstein injective over $A
\ltimes C$, and where
\[
  \id_A(C\LTensor_A Y) \leqslant \CGid M.
\]
\end{Lemma}

\begin{Lemma} 
\label{lem:pd_finite}
Let $C$ be a semi-dualizing module for $A$, let $M$ be a complex of
$A$-modules with homology bounded to the right and $\pd_A M < \infty$,
and let $H$ be an $A$-module which is Gorenstein injective over $A
\ltimes C$.  Then
\[
  \H_{-(j+1)}\RHom_A(M,H) \,=\, 0
\]
for $j \geqslant \sup \{\, i \,|\, \H_i\!M \not= 0 \,\}$. 
\end{Lemma} 

\begin{proof}  
Since $M$ has homology bounded to the right and $\pd_A M < \infty$,
there exists a bounded projective resolution $P$ of $M$, and
\[
  \H_{-(j+1)}\RHom_A(M,H) \,\cong\, \Ext^1_A(C^P_j,H)
\]
where $C^P_j$ is the $j$'th cokernel of $P$.  Since
\[
  \cdots \longrightarrow P_{j+1} \longrightarrow P_j 
  \longrightarrow C^P_j \rightarrow 0
\]
is a projective resolution of $C^P_j$ and since $P$ is bounded, we have
$\pd_A C^P_j < \infty$. Hence it is enough to show
\[
  \Ext^1_A(M,H) = 0
\]
for each $A$-module $M$ with $\pd_A M < \infty$.  

To prove this, we first argue that if $I$ is any injective
$A$-module then
\[
  \Ext^i_A(M,\Hom_A(C,I))=0
\]
for $i > 0$.  For this, note that we have
\begin{align*}
  \RHom_A(M,\Hom_A(C,I)) 
  & \cong \RHom_A(M,\RHom_A(C,I)) \\
  & \stackrel{\rm (a)}{\cong} \RHom_A(M \LTensor_A C,I) \\
  & \stackrel{\rm (b)}{\cong} \RHom_A(C,\RHom_A(M,I)) \\
  & \cong \RHom_A(C,\Hom_A(M,I))
\end{align*}
where (a) and (b) are by adjunction, and consequently,
\begin{equation}
\label{equ:Henrikdagger1}
  \Ext^i_A(M,\Hom_A(C,I)) \cong \Ext^i_A(C,\Hom_A(M,I))
\end{equation}
for each $i$.  The condition $\pd_A M < \infty$ implies $\id_A
\Hom_A(M,I) < \infty$, and therefore $\Hom_A(M,I)$ belongs to
${}_{C}\cB(A)$ by \cite[prop.\ (4.4)]{Wintherpaper}.  Thus
\cite[obs.\ (4.10)]{Wintherpaper} implies that the right hand side of
\eqref{equ:Henrikdagger1} is zero for $i > 0$.

Now set $n = \pd_A M$. Repeated use of lemma \ref{dualt_PJ-Lemma10}
shows that there is an exact sequence of $A$-modules
\begin{equation}
\label{equ:Henrikddagger1}
  0 \rightarrow H^{\prime} \longrightarrow \Hom_A(C,I_{n-1}) 
  \longrightarrow \cdots \longrightarrow \Hom_A(C,I_0)
  \longrightarrow H \rightarrow 0, 
\end{equation}
where $I_0, \ldots, I_{n-1}$ are injective $A$-modules. Applying
$\Hom_A(M,-)$ to \eqref{equ:Henrikddagger1} and using
$\Ext^i_A(M,\Hom_A(C,I_j)) = 0$ for each $i > 0$ and each $j$,
we obtain
\[
  \Ext^1_A(M,H) \,\cong\, \Ext^{n+1}_A(M,H^{\prime}) = 0
\]
as desired.  Here the last equality holds because $\pd_A M = n$.
\end{proof}  

\begin{Lemma}
\label{lem:id_vs_CGid}
Let $C$ be a semi-dualizing module for $A$ and let $M$ be in
${}_{C}\cA(A)$.  Set $s = \sup \{\, i \,|\, \H_i\!M \not= 0 \,\}$ and
suppose that $M$ satisfies
\[
  \H_{-(s + 1)} \RHom_A(M,H) = 0
\]
for each $A$-module $H$ which is Gorenstein injective over $A \ltimes
C$.  Then
\[
  \CGid M  = \id_A (C\LTensor_AM).
\]
\end{Lemma}

\begin{proof}
To prove the lemma's equation, let us first prove the inequality
$\leqslant$.  Let $t = \sup \{\, i \,|\, \H_i(C \LTensor_A M) \not= 0
\,\}$ and $n = \id_A(C\LTensor_AM)$.  We may clearly suppose that $n$
is finite.  Let
\[
  J = \cdots \longrightarrow 0 \longrightarrow J_t
  \longrightarrow \cdots \longrightarrow \cdots \longrightarrow
  J_{-n} \longrightarrow 0 \longrightarrow \cdots
\]
be an injective resolution of $C\LTensor_AM$.  The complex $M$ is in
${}_{C}\cA(A)$, so we get the first $\cong$ in
\[
  M \cong \RHom_A(C,C \LTensor_A M) \cong \Hom_A(C,J).
\]
Lemma \ref{lem:A_and_C_Gproj}(2) implies that 
$\Hom_A(C,J)$ is a complex of Gorenstein injective modules
over $A \ltimes C$.  Since $\Hom_A(C,J)_{\ell} = \Hom_A(C,J_{\ell}) =
0$ for $\ell < -n$, we see that
\[
  \CGid M \leqslant n = \id_A(C\LTensor_AM). 
\]

Let us next prove the inequality $\geqslant$.  Recall that $s = \sup
\{ i \,|\, \H_i\!M \not= 0 \,\}$.  We may clearly suppose that
$\CGid M$ is finite.  By lemma \ref{lem:triangle} there is a
distinguished triangle in $\D(A)$,
\begin{equation}
\label{equ:Henrikast}
  \Sigma^s H 
  \longrightarrow Y 
  \stackrel{f}{\longrightarrow} M 
  \longrightarrow,
\end{equation}
where $H$ is an $A$-module which is Gorenstein injective over $A
\ltimes C$, and where
\begin{equation}
\label{equ:Henrikastast}
  \id_A(C\LTensor_A Y) \leqslant \CGid M.
\end{equation}

Applying $\RHom_A(M,-)$ to \eqref{equ:Henrikast} gives another
distinguished triangle whose long exact homology sequence contains
\[
  \H_0\RHom_A(M,Y) \longrightarrow \H_0\RHom_A(M,M) \longrightarrow
  \H_{-1}\RHom_A(M,\Sigma^s H)
\]
which can also be written
\[
  \Hom_{\Dsmall(A)}(M,Y) \longrightarrow \Hom_{\Dsmall(A)}(M,M) 
  \longrightarrow \H_{-(s + 1)}\RHom_A(M,H) = 0,
\]
where the last zero comes from the assumptions on $M$. Consequently,
there exists a morphism $g : M \longrightarrow Y$ in $\D(A)$ with 
$fg = 1_M$.  That is, the distinguished triangle \eqref{equ:Henrikast}
is split, so $Y \cong \Sigma^s H \oplus M$.

This implies
\[
  C \LTensor_A Y \cong (C \LTensor_A \Sigma^s H)
                       \oplus (C \LTensor_A M)
\]
from which clearly follows
\begin{equation}
\label{equ:Henrikastastast}
  \id_A(C \LTensor_A M) \leqslant \id_A(C \LTensor_A Y).
\end{equation}
Combining the inequalities \eqref{equ:Henrikastast} and
\eqref{equ:Henrikastastast} now shows
\[
  \id_A(C \LTensor_A M) \leqslant \CGid M
\]
as desired.
\end{proof}

\begin{Proposition}
\label{pro:id_inequality}
Assume that the ring $A$ is local.  Let $C$ be a semi-dualizing module
for $A$ and let $M$ be a complex of $A$-modules with homology bounded
to the right and $\fd_A M < \infty$. Then
\begin{align*} 
  \id_A C \leqslant \CGid M + \width_A M.
\end{align*}
\end{Proposition}

\begin{proof}
Denote by $k$ the residue class field of $A$.  Since $\fd_A M <
\infty$, the isomorphism \cite[(A.4.23)]{Wintherbook} gives
\[
  \RHom_A(k,C \LTensor_A M) \cong \RHom_A(k,C) \LTensor_A M.
\]
This implies (a) in
\begin{align*} 
  \lefteqn{
           \inf \{\, i \,|\, \H_i \RHom_A(k,C \LTensor_A M) \not= 0 \,\}
          } \;\;\;\;\; 
  & \\
  & \stackrel{\rm (a)}{=} 
    \inf \{\, i \,|\, \H_i(\RHom_A(k,C) \LTensor_A M) \not= 0 \,\} \\
  & \stackrel{\rm (b)}{=}
    \inf \{\, i \,|\, \H_i \RHom_A(k,C) \not= 0 \,\}
    + \inf \{\, i \,|\, \H_i(M \LTensor_A k) \not= 0 \,\} \\
  & = - \id_A C + \width_A M,
\end{align*}
where (b) is by \cite[(A.7.9.2)]{Wintherbook}.  Consequently,
\begin{align*} 
  \id_A C 
  & = - \inf \{\, i \,|\, \H_i \RHom_A(k,C \LTensor_A M) \not= 0 \,\}
      + \width_A M \\
  & \leqslant \id_A(C \LTensor_A M) + \width_A M \\
  & = \CGid M + \width_A M.
\end{align*}
The last $=$ follows from lemmas \ref{lem:pd_finite} and
\ref{lem:id_vs_CGid}.  Lemma \ref{lem:pd_finite} applies to $M$
because $\fd_A M < \infty$ implies $\pd_A M < \infty$ by \cite[Seconde
partie, cor.\ (3.2.7)]{RaynaudGruson}, and lemma \ref{lem:id_vs_CGid}
applies because $\fd_A M < \infty$ implies $M \in {}_{C}\cA(A)$ by
\cite[prop.\ (4.4)]{Wintherpaper}.
\end{proof}

\begin{Lemma}
\label{lem:CGid_of_dual}
Let $C$ be a semi-dualizing module for $A$, let $I$ be a faithfully
injective $A$-module, and let $M$ be a complex of $A$-modules with
right-bounded homology.  Then
\[
  \CGid \Hom_A(M,I) = \CGfd M.
\]
\end{Lemma}

\begin{proof}
From lemma \ref{lem:induced_injectives}(1) follows that $E = \Hom_A(A
\skewtimes C,I)$ is a faithfully injective $(A \skewtimes C)$-module.
Hence
\[
  \Gid_{A \skewtimes C} \Hom_{A \skewtimes C}(M,E)
  = \Gfd_{A \skewtimes C} M
\]
follows from \cite[thm.\ (6.4.2)]{Wintherbook}.

But equation \eqref{equ:adjunction} in the proof of lemma
\ref{lem:induced_injectives} shows $\Hom_{A \skewtimes C}(M,E) \cong
\Hom_A(M,I)$, so accordingly,
\[
  \Gid_{A \skewtimes C} \Hom_A(M,I) = \Gfd_{A \skewtimes C} M.
\]
\end{proof}

\begin{Proposition} 
\label{pro:fd_inequality}
Assume that the ring $A$ is local.  Let $C$ be a semi-dualizing
module for $A$ and let $N$ be a complex of $A$-modules with homology
bounded to the left and $\id_A N < \infty$. Then
\begin{align*} 
  \id_A C \leqslant \CGfd N + \textnormal{depth}_A N.
\end{align*}
\end{Proposition} 

\begin{proof}
Apply Matlis duality and lemma \ref{lem:CGid_of_dual} to proposition
\ref{pro:id_inequality}.
\end{proof}

\section{Properties of the Cohen-Macaulay dimensions}
\label{sec:CMdims}

\begin{Theorem}
\label{thm:CMdims_finite}
Assume that the ring $A$ is local with residue class field $k$.  Then
the following conditions are equivalent.
\begin{enumerate}

  \item  $A$ is a Cohen-Macaulay ring with a dualizing module.

\smallskip

  \item  $\CMid_A M < \infty$ holds when $M$ is any complex of
         $A$-modules with bounded homology.

\smallskip

  \item  There is a complex $M$ of $A$-modules with bounded homology,
         $\CMid_A M < \infty$, $\fd_A M < \infty$, and $\width_A M <
         \infty$. 

\smallskip

  \item  $\CMid_A k < \infty$.

\smallskip

  \item  $\CMpd_A M < \infty$ holds when $M$ is any complex of
         $A$-modules with bounded homology.

\smallskip

  \item  There is a complex $M$ of $A$-modules with bounded homology,
         $\CMpd_A M < \infty$, $\id_A M < \infty$, and $\depth_A M <
         \infty$. 

\smallskip

  \item  $\CMpd_A k < \infty$.

\smallskip

  \item  $\CMfd_A M < \infty$ holds when $M$ is any complex of
         $A$-modules with bounded homology.

\smallskip

  \item  There is a complex $M$ of $A$-modules with bounded homology,
         $\CMfd_A M < \infty$, $\id_A M < \infty$, and $\depth_A M <
         \infty$. 

\smallskip

  \item  $\CMfd_A k < \infty$.

\smallskip

\end{enumerate}
\end{Theorem}

\begin{proof}
Let us prove that conditions (1), (2), (3), and (4) are equivalent.
Similar proofs give that so are (1), (5), (6), and (7) as well as (1),
(8), (9), and (10).

\medskip
\noindent
(1) $\Rightarrow$ (2).  Let $A$ be Cohen-Macaulay with dualizing
module $C$.  Then $A \skewtimes C$ is Gorenstein by lemma
\ref{lem:AxC_Gorenstein}.  If $M$ is a complex of $A$-modules with
bounded homology, then $M$ is also a complex of $(A \skewtimes
C)$-modules with bounded homology, so
\[
  \Gid_{A \skewtimes C} M < \infty
\]
by \cite[thm.\ (6.2.7)]{Wintherbook}.  But as $C$ is in particular a
semi-dualizing module, the definition of $\CMid$ then implies
\[
  \CMid_A M < \infty.
\]

\medskip
\noindent
(2) $\Rightarrow$ (3) and (2) $\Rightarrow$ (4).  Trivial. 

\medskip
\noindent
(3) $\Rightarrow$ (1).  For $\CMid_A M < \infty$, the definition of
$\CMid$ implies that $A$ has a semi-dualizing module $C$ with 
\[
  \CGid M < \infty.  
\]
But when 
\[
  \fd_A M < \infty \; \; \mbox{and} \; \; \width_A M < \infty
\]
also hold, then proposition \ref{pro:id_inequality} implies 
\[
  \id_A C < \infty.
\]
So $A$ is Cohen-Macaulay with dualizing module $C$.

\medskip
\noindent
(4) $\Rightarrow$ (1).  When $\CMid_A k < \infty$ then $A$ has a
semi-dualizing module $C$ with 
\[
  \CGid k < \infty.  
\]
Then $\RHom_{A \skewtimes C}(E_{A \skewtimes C}(k),k)$ has bounded
homology by \cite[thm.\ 2.22]{HHGorensteinHomDim}.  However,
\begin{align*}
  \RHom_{A \skewtimes C}(E_{A \skewtimes C}(k),k)
  & \stackrel{\rm (a)}{\cong} 
    \RHom_{A \skewtimes C}(k^{\vee},E_{A \skewtimes C}(k)^{\vee}) \\
  & \cong
    \RHom_{A \skewtimes C}(k,\widehat{A \skewtimes C}) \\
  & \stackrel{\rm (b)}{\cong} 
    \RHom_{A \skewtimes C}(k,A \skewtimes C) 
      \otimes_{A \skewtimes C} \widehat{A \skewtimes C}, \\
\end{align*}
where (a) is by Matlis duality and (b) is by
\cite[(A.4.23)]{Wintherbook}.  Since $\widehat{A \skewtimes C}$ is
faithfully flat over $A \skewtimes C$, it follows that $\RHom_{A
\skewtimes C}(k,A \skewtimes C)$ also has bounded homology, whence $A
\skewtimes C$ is a Gorenstein ring.  But then $A$ is a Cohen-Macaulay
ring with dualizing module $C$ by lemma \ref{lem:AxC_Gorenstein}.
\end{proof}

\begin{Remark}
\label{rmk:CMdims_finite}
In condition (3) of the above theorem, one could consider for $M$
either the ring $A$ itself, or the Koszul complex $K(x_1, \ldots,
x_r)$ on any sequence of elements $x_1, \ldots, x_r$ in the maximal
ideal.  These complexes satisfy $\fd_A M < \infty$ and $\width_A M <
\infty$, and so either of the conditions
\[
  \CMid_A A < \infty 
  \; \; \mbox{and} \; \; \CMid_A K(x_1, \ldots, x_r) < \infty
\]
is equivalent to $A$ being a Cohen-Macaulay ring with a dualizing
module.

Similarly, in conditions (6) and (9), one could consider for $M$
either the injective hull of the residue class field, $E_A(k)$, or a
dualizing complex $D$ (if one is known to exist).  These complexes
satisfy  $\id_A M < \infty$ and $\depth_A M < \infty$, and so either
of the conditions 
\[
  \CMpd_A E_A(k) < \infty
  \; \; \mbox{and} \; \; \CMpd_A D < \infty
\]
and
\[
  \CMfd_A E_A(k) < \infty
  \; \; \mbox{and} \; \; \CMfd_A D < \infty
\]
is equivalent to $A$ being a Cohen-Macaulay ring with a dualizing
module.
\end{Remark}

The following results use $\CMdim$, the Cohen-Macaulay dimension
introduced by Gerko in \cite{Gerko}, and $\Gdim$, the G-dimension
originally introduced by Auslander and Bridger in
\cite{AuslanderBridger}. 

\begin{Lemma}
\label{lem:Gerko_inequality}
Assume that the ring $A$ is local.  Let $C$ be a semi-dualizing module
for $A$ and let $M$ be a finitely generated $A$-module.  If
\[
  \CGpd M < \infty
\]
then
\[
  \CMdim_A M = \CGpd M.
\]
\end{Lemma}

\begin{proof}
Combining \cite[proof of thm.\ 3.7]{Gerko} with \cite[def.\
3.2']{Gerko} shows 
\[
  \CMdim_A M \leqslant \CGpd M.
\]
So $\CGpd M < \infty$ implies $\CMdim_A M < \infty$ and hence
\[
  \CMdim_A M = \depth_A A - \depth_A M
\]
by \cite[thm.\ 3.8]{Gerko}.  

On the other hand,
\[
  \CGpd M = \CGdim_A M
\]
by \cite[prop.\ 3.1]{HHPJ1}, where $\CGdim_A$ is the dimension
introduced in \cite[def.\ (3.11)]{Wintherpaper} under the name
$G$-$\dim_C$.  So $\CGdim_A M$ is finite and hence
\[
  \CGdim_A M = \depth_A A - \depth_A M
\]
by \cite[thm.\ (3.14)]{Wintherpaper}.

Combining the last three equations shows
\[
  \CMdim_A M = \CGpd M
\]
as desired.
\end{proof}

%
%
%
%

\begin{Theorem}
\label{thm:CMpd_and_CMdim}
Assume that the ring $A$ is local and let $M$ be a finitely generated
$A$-module.  Then
\[
  \CMdim_A M \leqslant \CMpd_A M \leqslant \Gdim_A M,
\]
and if one of these numbers is finite then the inequalities to its
left are equalities.
\end{Theorem}

\begin{proof}
The first inequality is clear from lemma \ref{lem:Gerko_inequality},
since $\CMpd_A M$ is defined as the infimum of all $\CGpd M$.

For the second inequality, note that the ring $A$ is itself a
semi-dualizing module, so the definition of $\CMpd$ gives
$\leqslant$ in
\[
  \CMpd_A M 
  \leqslant \Gpd_{A \skewtimes A} M
  = \Gpd_A M
  = \Gdim_A M,
\]
where the first $=$ is by \cite[cor.\ 2.17]{HHPJ1} and the second $=$
holds because $M$ is finitely generated.

Equalities: If $\Gdim_A M < \infty$ then $\CMdim_A M < \infty$ by
\cite[thm.\ 3.7]{Gerko}.  But $\Gdim_A M < \infty$ implies
\[
  \Gdim_A M = \depth_A A - \depth_A M
\]
by \cite[thm.\ (2.3.13)]{Wintherbook}, and similarly, $\CMdim_A M <
\infty$ implies 
\[
  \CMdim_A M = \depth_A A - \depth_A M
\]
by \cite[thm.\ 3.8]{Gerko}.  So it follows that $\CMdim_A M = \Gdim_A
M$, and hence both inequalities in the theorem must be equalities.

If $\CMpd_A M < \infty$ then by the definition of $\CMpd$ there exists
a semi-dualizing module $C$ over $A$ with $\CGpd M < \infty$.  But
by lemma \ref{lem:Gerko_inequality}, any such $C$ has
\[
  \CMdim_A M = \CGpd M,
\]
and so it follows that the first inequality in the theorem is an
equality. 
\end{proof}

Since much is known about $\CMdim_A$, this theorem has several
immediate consequences for $\CMpd_A$.  The following is even clear
from the proof of the theorem.

\begin{Theorem}
[Auslander-Buchsbaum formula]
\label{thm:AB}
Assume that the ring $A$ is local and let $M$ be a finitely generated
$A$-module.  If $\CMpd_A M < \infty$, then
\[
  \CMpd_A M = \depth_A A - \depth_A M.
\]
\end{Theorem}

\bigskip

\noindent
{\bf Acknowledgement.}  
The diagrams were typeset with Paul Taylor's {\tt diagrams.tex}.

\end{document}